\documentclass{amsart}
\usepackage{bbm}
\usepackage{mathrsfs}
\usepackage{amsmath}
\usepackage{paralist}
\usepackage{graphics}
\usepackage{epsfig}

\usepackage{amsfonts}
\usepackage{CJK}
\usepackage{fancyhdr}
\usepackage{graphicx}
\usepackage[dvips,usenames]{color}
\usepackage{titletoc}
\usepackage{latexsym}
\usepackage{amssymb}
\usepackage{multicol}
\usepackage{graphics}
\usepackage{subfigure}
\usepackage{indentfirst}
\usepackage{cases}
\usepackage{curves}

\newtheorem{theorem}{Theorem}[section]
\newtheorem{corollary}[theorem]{Corollary}
\newtheorem{lemma}[theorem]{Lemma}
\newtheorem{definition}[theorem]{Definition}
\newtheorem{proposition}{Proposition}

\newtheorem{example}[theorem]{Example}
\numberwithin{equation}{section}

\def\x#1{(\ref{#1})}
\def\R{{\Bbb R}}

\def\Z{{\Bbb Z}}

\newcommand{\cF}{\mathcal{F}}

\newcommand{\cS}{\mathcal{S}}
\newcommand{\cU}{\mathcal{U}}
\newcommand{\cW}{\mathcal{W}}
\newcommand{\cI}{\mathcal{I}}

\newcommand{\ned}{~nonuniform exponential dichotomy}

\newcommand{\nds}{nonuniform dichotomy spectrum~}

\def\bc{\begin{center}}
\def\ec{\end{center}}
\def\ba{\begin{array}}
\def\ea{\end{array}}
\def\be{\begin{equation}}
\def\ee{\end{equation}}
\def\bea{\begin{eqnarray}}
\def\eea{\end{eqnarray}}
\def\beaa{\begin{eqnarray*}}
\def\eeaa{\end{eqnarray*}}

\def\ben{\begin{enumerate}}
\def\een{\end{enumerate}}
\def\hh{\!\!\!\!}

\def\EQ{\hh & = & \hh}

\def\LE{\hh & \le & \hh}

\def\nn{\nonumber}
\def\oo{\infty}
\def\ifl{\iffalse}
\def\lb{\label}
\def\prf{\mbox{\bf Proof.~}}

\DeclareMathOperator{\im}{im}

\title[]{Nonuniform dichotomy spectrum and reducibility for nonautonomous difference equations}

\author[Jifeng Chu,\quad Hailong Zhu, \quad Stefan Siegmund \quad and \quad Yonghui Xia]
{Jifeng Chu$^{1,2}$, \quad Hailong Zhu$^1$, \quad Stefan
Siegmund$^2$,\quad Yonghui Xia$^3$}
\address{$^1$ Department of Mathematics, College of Science,
Hohai University, Nanjing 210098, China}
\address{$^2$ Center for Dynamics \& Institute for Analysis, Department of Mathematics, TU Dresden, Germany}
\address{$^3$ College of Mathematics, Physics and Information Engineering, Zhejiang
Normal University, Jinhua, China}
\email{jifengchu@126.com (J.\ Chu)}
\email{hai-long-zhu@163.com (H. Zhu)}
\email{stefan.siegmund@tu-dresden.de (S.\ Siegmund)}
\email{yhxia@zjnu.cn; xiadoc@163.com (Y.\ Xia)}
\thanks{Jifeng Chu was supported by the National Natural Science
Foundation of China (Grant No. 11171090, No. 11271078 and No.
11271333), China Postdoctoral Science Foundation funded project
(Grant No.2012T50431) and the Alexander von Humboldt Foundation of
Germany. Hailong Zhu was supported by the National Natural Science
Foundation of China (Grant No. 11301001), NSF of the Educational
Bureau of Anhui province (Grant NO. 10040606Q01, NO. 1208085QA11,
NO. 1208085QG131) }

\subjclass[2000]{37D25; 37B55} \keywords{Dichotomy spectrum;
nonuniform exponential dichotomy; reducibility.}

\begin{document}

\begin{abstract} For nonautonomous linear difference equations, we
introduce the notion of the so-called nonuniform dichotomy spectrum
and prove a spectral theorem. Moreover, we introduce the notion of weak
kinematical similarity and prove a reducibility result by the spectral theorem.
\end{abstract}

\maketitle

\section{\bf Introduction}
\setcounter{equation}{0} \noindent Let $A_{k}\in \R^{N \times
N},~k\in \Z$, be a sequence of invertible matrices. In this
paper, we consider the following nonautonomous linear difference
equations \be\lb{le}x_{k+1}=A_{k}x_{k},\ee where $x_{k}\in \R^{N},k\in\Z$.  Let $\Phi:\Z\times\Z\rightarrow\R^{N\times N},
(k,l)\mapsto\Phi(k,l)$, denote the evolution operator of \x{le},
i.e.,
\[\Phi(k,l)=\begin{cases}
A_{k-1}\cdots  A_{l},  &{\rm for}\quad k>l,\\
 {\rm Id}, &{\rm for}\quad k=l,\\
 A^{-1}_{k}\cdots  A^{-1}_{l-1}, &{\rm for}\quad k<l.
\end{cases}\]Obviously, $\Phi(k,m)\Phi(m,l)=\Phi(k,l)$, $k,m,l\in \Z$, and
$\Phi(\cdot,l)\xi$ solves the initial value problem \x{le},
$x(l)=\xi$, for $l \in \Z$, $\xi\in\R^{N}$.

An invariant projector of \x{le} is defined to be a function
$P:\Z\rightarrow\R^{N\times N}$ of projections $P_{k}, k\in \Z$,
such that for each $P_k$ the following property holds
\[P_{k+1}A_{k}=A_{k}P_{k},\quad k\in \Z.\]
We say that {\rm\x{le}} admits an {\it exponential
dichotomy} if there exist an invariant projector $P$
and constants $0<\alpha<1,K\geq1$ such that
\be\lb{pl}\|\Phi(k,l)P_{l}\| \leq K\alpha^{k-l}, \quad k\geq l,\ee
and \be\lb{ql}\|\Phi(k,l)Q_l\|\leq K(\tfrac{1}{\alpha})^{k-l},\quad
k\leq l,\ee where $Q_l={\rm Id}-P_l$ is the complementary
projection.

The notion of exponential dichotomy was introduced by
Perron in \cite{per} and has attracted a lot of interest during the
last few decades because it plays an important role in the study of hyperbolic
dynamical behavior of differential equations and difference
equations. For example, see \cite{aul, gkk, ss} and the
references therein. We also refer to the books \cite{clr, cop, ms}
for details and further references related to exponential
dichotomies. On the other hand, during the last decade,
inspired both by the classical notion of exponential dichotomy and
by the notion of nonuniformly hyperbolic trajectory introduced by
Pesin (see \cite{bp02}), Barreira and Valls have
introduced the notion of nonuniform exponential dichotomies and have
developed the corresponding theory in a systematic way \cite{bv08,
bv-jde-09, bv-jde-10, bv-non, bv-dcds-12, bv-bsm, bv-dcds-12-2,
bv-dcds-13}. As explained by Barreira and Valls, in comparison to
the notion of exponential dichotomies, nonuniform exponential
dichotomy is a useful and weaker notion. A very general type of nonuniform
exponential dichotomy has been considered in \cite{bcv-2,bcv-1, chu}.

We say that {\rm\x{le}} admits a {\it nonuniform
exponential dichotomy} if there exist an invariant projector $P$ and
constants $0<\alpha<1,K\geq1,\varepsilon\geq 1$, such that
\be\lb{npl}\|\Phi(k,l)P_{l}\| \leq K\alpha^{k-l}\varepsilon^{l},
\quad k\geq l,\ee and \be\lb{nql}\|\Phi(k,l)Q_l\|\leq
K(\tfrac{1}{\alpha})^{k-l}\varepsilon^{l},\quad k\leq l.\ee
When $\varepsilon=1$, \x{npl}-\x{nql} become \x{pl}-\x{ql}, and
therefore a nonuniform exponential dichotomy
becomes an exponential dichotomy.
For example, given $\omega>a>0$, then the linear
equation
\be\lb{4.1}u_{k+1}=e^{-\omega+ak(-1)^{k}-a(k-1)(-1)^{(k-1)}}u_{k},\,\,
\,v_{k+1}=e^{\omega-ak(-1)^{k}+a(k-1)(-1)^{(k-1)}}v_{k}\ee admits a
nonuniform exponential dichotomy, but does not admit an exponential
dichotomy. In fact, we have
\[\Phi(k,l)P_{l}=\left(\begin{array}{lll} e^{-\omega(k-l-1)-a(k-l-1)(-1)^{k-1}-al(-1)^{(k-1)}+al(-1)^{l}} & 0
 \\
 0 &  0
 \end{array}\right)\]
with $P_{l}=\left(\begin{array}{lll} 1 & 0
 \\
 0 &  0
 \end{array}
 \right)$. Therefore \x{npl} holds with
\[K=e^{\omega-a}>1,\quad \alpha=e^{(-\omega+a)}\in(0,1),\quad \varepsilon=e^{2a}>1.\] Analogous arguments applied to
the second equation yield the estimate \x{nql}. Moreover, when
both $k$ and $l$ are even, we obtain the equality
\[\|\Phi(k,l)P_{l}\| = K\alpha^{k-l}\varepsilon^{l}, \quad k\geq
l,\] which means that the nonuniform part $\varepsilon^{l}=e^{2al}$
cannot be removed.

Although the notion of nonuniform exponential dichotomy has been studied in a
very wide range and many rich results have been obtained,
up to now there are no results on the spectral theory of
\x{le} in the setting of nonuniform exponential dichotomies. In this paper,
we establish the spectral theory in the setting of strong
nonuniform exponential dichotomies. We say that \x{le} admits a {\it strong nonuniform exponential
dichotomy} if it admits a \ned~with $\alpha\varepsilon^2<1$ in \x{npl}-\x{nql}. For example, if $\omega>5a$, then \x{4.1} admits a
strong nonuniform exponential dichotomy. We remark that the phrase ``strong
nonuniform exponential dichotomy" has been used in \cite{bv08}, however here we use
this notion in a different sense. Moreover, \cite[Theorem 1.4.2]{bp02} indicates that
the condition $\alpha\varepsilon^2<1$ is reasonable, which means that the constant $\varepsilon$
belongs to the interval $[1,\sqrt{1/\alpha})$.

Among the different
topics on classical exponential dichotomies, the dichotomy spectrum is very important and many results
have been obtained. We refer the reader to \cite{amz, am, cl94, ns, pot09, pot12, ss78, s-jdde, s-jlms} and the references therein.
The definition and investigation for finite-time hyperbolicity has also been
studied in \cite{bds, dkns, dps}.

This paper is organized as follows. In Section 2 we propose a definition of
spectrum based on strong nonuniform exponential dichotomies, which is called
nonuniform dichotomy spectrum. Such a spectrum can be seen as a generalization of Sacker-Sell spectrum.
We prove a nonuniform dichotomy spectral theorem. In Section 3 we prove a
reducibility result for \x{le} using the spectral result. Recall that system (\ref{le})
is reducible if it is kinematically similar to a block diagonal
system with blocks of dimension less than $N$.

\section{Nonuniform dichotomy spectrum}
\setcounter{equation}{0}

Consider the weighted system \be\lb{2.1} x_{k+1} =
\tfrac{1}{\gamma}A_{k}x_{k},\ee where $\gamma \in\R^+=(0,\oo)$. One
can easily see that
\[\Phi_\gamma(k,l) := (\tfrac{1}{\gamma})^{k-l}\Phi(k,l)\] is its evolution
operator. If for some $\gamma\in\R^+$, \x{2.1} admits a
nonuniform exponential dichotomy with projector $P_k$ and constants $K\geq 1,0<\alpha<1$ and
$\varepsilon\geq 1$, then $P_k$ is also invariant for \x{le}, that is
\[P_{k+1}A_{k}=A_{k}P_{k},\quad k\in \Z,\] and the
dichotomy estimates of \x{2.1} are equivalent to
\be\lb{npe}
\|\Phi(k,l)P_{l}\| \leq K(\gamma\alpha)^
{k-l}\varepsilon^{l}, \quad k\geq l,\ee
and
\be\lb{nqe}\|\Phi(k,l)Q_l\|\leq K(\gamma\tfrac{1}
{\alpha})^{k-l}\varepsilon^{l},\quad k\leq l.\ee

\begin{definition} \label{Def2.1}
The nonuniform dichotomy spectrum of {\rm(\ref{le})} is the set
\[\Sigma_{NED}(A) = \{ \gamma \in\R^+:\, \x{2.1} \mbox{ admits no strong nonuniform exponential dichotomy} \},\]
and the \emph{resolvent set} $\rho_{NED}(A) =\R^+\setminus
\Sigma_{NED}(A)$ is its complement. The dichotomy spectrum
of {\rm(\ref{le})} is the set
\[\Sigma_{ED}(A) = \{ \gamma \in\R^+:\, \x{2.1} \mbox{ admits no exponential dichotomy} \},\]
and $\rho_{ED}(A) =\R^+\setminus \Sigma_{ED}(A)$.
\end{definition}

\begin{proposition} \label{pro2.1}
$\Sigma_{NED}(A)\subset \Sigma_{ED}(A)$.
\end{proposition}

\noindent{\bf Proof.} For each $\gamma \in\rho_{ED}(A) $, the
weighted system (\ref{2.1}) admits an exponential dichotomy. Consequently, the weighted
system (\ref{2.1}) admits a strong \ned. Thus, $\gamma \in\rho_{NED}(A)$,
which implies that $\rho_{ED}(A)\subset\rho_{NED}(A) $, and therefore
$\Sigma_{NED}(A)\subset \Sigma_{ED}(A)$.   $\Box$

Let us define for $\gamma \in\rho_{NED}(A)$
\[\cS_\gamma := \{(l,\xi) \in \Z \times \R^N:
\sup_{k\geq l}\|\Phi(k,l)\xi\|\gamma^{-k}\varepsilon^{-l}<\oo\},\]
and
\[\cU_\gamma:= \{(l,\xi) \in \Z \times \R^N:
\sup_{k\leq l}\|\Phi(k,l)\xi\|\gamma^{-k}\varepsilon^{-l}<\oo\},\] where
$\varepsilon$ is the constant in \x{npe}-\x{nqe}.
 One may readily verify that
$\cS_{\gamma}$ and $\cU_{\gamma}$ are invariant vector bundles of
\x{le}, here we say that a nonempty set $\mathcal{W} \subset \Z
\times \R^N$ is an invariant vector bundle of \x{le} if (a) it is
{\rm invariant}, i.e., $(l,\xi) \in \mathcal{W} \;\Rightarrow\;
(k,\Phi(k,l)\xi) \in \mathcal{W}$ for all $k \in \Z$; and (b) for
every $l \in \Z$ the {\rm fiber} $\mathcal{W}(l) = \{ \xi \in \R^N
\,:\, (l,\xi) \in \mathcal{W} \}$ is a linear subspace of $\R^N$.

As a first glance, $\cS_\gamma$ and $\cU_\gamma$ are not well defined
because they seem to depend on the constant $\varepsilon$, which may
be not unique in \x{npe}-\x{nqe}. However, the following
result ensures that $\cS_\gamma$ and $\cU_\gamma$ are well defined and they
do not depend on the choice of the constant $\varepsilon$. First we recall that
the invariant projector $P$ is unique for \x{le} and \x{2.1} following the arguments in \cite[Chapter 2]{cop}.
Although the arguments in \cite{cop} are done in the setting of exponential dichotomies,
it is not difficult to verify that they are also applicable to the case of nonuniform exponential dichotomies.

\begin{lemma}  \label{lem2.1}
Assume that {\rm \x{2.1}} admits a strong \ned ~with invariant projector $P$ for
$\gamma \in\R^+.$ Then
\[\cS_\gamma = \im P
\;,\quad
\cU_\gamma = \ker P
\quad \text{and} \quad
\cS_\gamma \oplus \cU_\gamma = \Z \times \R^N
\;. \]
\end{lemma}

\noindent\prf{We show only $\cS_\gamma = \im P$. The fact $\cU_\gamma = \ker P$ is analog and the fact
$\cS_\gamma \oplus \cU_\gamma = \Z \times \R^N$ is clear.

First we show $\cS_\gamma\subset\im P.$ Let $l\in\Z$ and $\xi\in\cS_\gamma(l).$ Then there exists
a positive constant $C$ such that
\[\|\Phi(k,l)\xi\|\leq C\gamma^{k}\varepsilon^{l},\quad k\geq l.\]
We write $\xi=\xi_1+\xi_2$ with $\xi_1\in {\rm im} P_l$ and
$\xi_2\in {\rm ker} P_l.$ We show that $\xi_2=0.$ The invariance of
$P$ implies for $k\in\Z$, we have the identity
\[\xi_2=\Phi_\gamma(l,k)\Phi_\gamma(k,l)Q_l\xi=\Phi_\gamma(l,k)Q_k\Phi_\gamma(k,l)\xi.\]
Since \x{2.1} admits a strong \ned, the following inequality holds
\[\|\Phi_\gamma(l,k)Q_k\|\leq K(\tfrac{1}{\alpha})^{l-k}\varepsilon^{k}.\]
Thus  \beaa\|\xi_2\|\LE
K(\tfrac{1}{\alpha})^{l-k}\varepsilon^{k}\|\Phi_\gamma(k,l)\xi\|\\
\EQ K(\alpha\varepsilon)^{k-l}\varepsilon^{l}(\tfrac{1}{\gamma})^{k-l}\|\Phi(k,l)\xi\|\\
\LE CK(\alpha\varepsilon)^{k-l}\varepsilon^{2l}(\tfrac{1}{\gamma})^{k-l}\gamma^k\\
\EQ CK(\alpha\varepsilon)^{k-l}\varepsilon^{2l}\gamma^l\quad
k\geq l,\eeaa which implies that $\xi_2=0$ by letting $k\rightarrow \oo$, since $\alpha\varepsilon<1$.

Next we show ${\rm \im} P\subset \cS_\gamma.$ Let $l\in\Z$ and
$\xi\in {\rm im} P_l,$ i.e., $P_l\xi=\xi.$ The \ned ~implies that
\[\|\Phi_\gamma(k,l)\xi\|\leq K\alpha^{k-l}\varepsilon^{l}\|\xi\|\leq K\varepsilon^{l}\|\xi\|,\quad k\geq l,\]
since $\alpha<1$, which implies that
\[\|\Phi(k,l)\xi\|\leq K\gamma^{k-l}\varepsilon^{l}\|\xi\|,\]
and hence $\xi\in\cS_\gamma(l).$} $\Box$

\begin{lemma}\label{lemma2.2}  The resolvent set is open, i.e., for every $\gamma \in \rho_{NED}(A),$ there
exists a constant $\beta =\beta(\gamma)\in(0,1)$ such that $(\beta\gamma,
\frac{1}{\beta}\gamma) \subset \rho_{NED}(A)$. Furthermore,
\[\cS_\zeta = \cS_\gamma
  \;\quad and \quad
  \cU_\zeta =  \cU_\gamma
  \quad \text{for} \quad
 \zeta \in(\beta\gamma,
\tfrac{1}{\beta}\gamma)
  \;.
  \]
 \end{lemma}

\noindent\prf{Let $\gamma \in \rho_{NED}(A)$. Then \x{2.1} admits a strong
\ned, i.e., the estimates \x{npe}-\x{nqe} hold with an invariant
projector $P$, constants $K \geq 0$, $0<\alpha <1$ and
$\varepsilon\geq 1$. For $\beta:= \sqrt{\alpha}\in(0,1)$ and $\zeta \in
(\beta\gamma, \frac{1}{\beta}\gamma)$ we have \[\Phi_\zeta(k,l) =
(\tfrac{\gamma}{\zeta})^{k-l}\Phi_\gamma(k,l).\] Now $P$ is also an invariant
projector for \[x_{k+1} = \tfrac{1}{\zeta}A_{k}x_{k}.\] Moreover, we have
the estimates
\[\|\Phi_\zeta(k,l)P_{l}\| \leq
 K(\tfrac{\gamma}{\zeta}\alpha)^{k-l}\varepsilon^{l} \leq K\beta^{k-l}\varepsilon^{l}, \quad k\geq l,\]
 and
\[\|\Phi_\zeta(k,l)Q_{l}\|
 \leq
K(\tfrac{\gamma}{\zeta}\tfrac{1}{\alpha})^{k-l}\varepsilon^{l} \leq
K(\tfrac{1}{\beta})^{k-l}\varepsilon^{l},\quad k\leq l.\]
Hence $\zeta \in \rho_{NED}(A)$. Therefore, $\rho_{NED}(A)$ is an
open set. Using Lemma \ref{lem2.1}, we know that $\cS_\zeta = \cS_\gamma$ and $\cU_\zeta =
\cU_\gamma$.} $\Box$

\begin{corollary}\label{cor2.1}  \,\,
 $\Sigma_{NED}(A)$ is a closed set.
 \end{corollary}

Using the facts proved above, we can obtain the following result, whose proof
is similar as \cite[Lemma 2.2]{as02}, and therefore we omit the proof here.

\begin{lemma}  \label{lemma2.3}
Let $\gamma_1,\gamma_2 \in \rho_{NED}(A)$ with $\gamma_1 <
\gamma_2$. Then $\cF
= \cU_{\gamma_1} \cap \cS_{\gamma_2}$ is an invariant vector bundle
which satisfies exactly one of the following two alternatives and
the statements given in each alternative are equivalent:

\begin{tabular}{ll}
\hspace{1.5cm}Alternative I & \hspace{1.5cm}Alternative II
\\[1ex]
{\rm(A)} $\cF = \Z \times \{0\}$. & {\rm(A')} $\cF \not= \Z \times
\{0\}$.
\\[0.5ex]
{\rm(B)} $[\gamma_1,\gamma_2] \subset \rho_{NED}(A)$. & {\rm(B')}
There is a $\zeta \in (\gamma_1,\gamma_2) \cap \Sigma_{NED}(A)$.
\\[0.5ex]
{\rm(C)} $\cS_{\gamma_1} = \cS_{\gamma_2}$ and $\cU_{\gamma_1} =
\cU_{\gamma_2}$. \hspace{6mm} & {\rm(C')} $\dim \cS_{\gamma_1} <
\dim \cS_{\gamma_2}$.
\\[0.5ex]
{\rm(D)} $\cS_\gamma = \cS_{\gamma_2}$ and $\cU_\gamma =
\cU_{\gamma_2}$ & {\rm(D')} $\dim \cU_{\gamma_1} > \dim
\cU_{\gamma_2}$.
\\
\hspace{6.5mm}for $\gamma \in [\gamma_1,\gamma_2]$. &
\end{tabular}
\end{lemma}

Now we are in a position to state and prove the nonuniform dichotomy spectral theorem which
will be essential to prove the reducibility result in Section 3. The proof follows
the idea and technique of the classical dichotomy spectrum proposed
in \cite{phd}, we present the details for the
reader's convenience.

\begin{theorem}\label{main21} The
\nds $\Sigma_{NED}(A)$ of {\rm\x{le}} is the disjoint union of $n$
closed intervals {\rm(}called {\it spectral intervals}{\rm)} where $0 \leq n
\leq N$, i.e., $\Sigma_{NED}(A) = \emptyset$ or $\Sigma_{NED}(A)
=\R^+$ or one of the four cases
 \begin{equation*}\label{spectrumNED}
\Sigma_{NED}(A) =
  \left\{
  \begin{matrix}
    [a_1,b_1]
    \\
    \text{or}
    \\
    (0,b_1]
  \end{matrix}
  \right\}
  \cup [a_2,b_2] \cup \cdots \cup [a_{n-1},b_{n-1}] \cup
  \left\{
  \begin{matrix}
    [a_n,b_n]
    \\
    \text{or}
    \\
    {[}a_n,\infty)
  \end{matrix}
  \right\},
  \end{equation*}
where $0<a_1 \leq b_1 < a_2 \leq b_2 < \cdots < a_n \leq b_n$.
Choose a \be\lb{bc1}\gamma_0 \in \rho_{NED}(A)~ with~
(0,\gamma_0)\subset \rho_{NED}(A)~ if~ possible, \ee otherwise
define $\cU_{\gamma_0}:=\Z \times \R^{N}$, $\cS_{\gamma_0}:=\Z
\times \{0\}$. Choose a \be\lb{bc2}\gamma_n \in \rho_{NED}(A)~ {
with}~ (\gamma_n,+\oo)\subset \rho_{NED}(A)~ if~ possible, \ee
otherwise define $\cU_{\gamma_n}:=\Z \times \{0\}$,
$\cS_{\gamma_0}:=\Z \times \R^{N}$. Then the sets
\[\cW_0=\cS_{\gamma_0}\quad{\rm and}\quad\cW_{n+1}=\cS_{\gamma_n}\]
are {invariant vector bundles} of {\rm \x{le}}. For $n\geq 2$,
choose $\gamma_i \in \rho_{NED}(A)$ with \be\lb{bc3}b_i < \gamma_i <
a_{i+1}\quad {for} \quad i=1,\ldots,n-1,\ee then for every
$i=1,\ldots,n-1$ the intersection \[\cW_i=\cU_{\gamma_{i-1}}\cap
\cS_{\gamma_{i}}\] is an invariant vector bundle of {\rm \x{le}}
with $\dim \cW_i \geq 1$. The invariant vector bundles $\cW_i,
i=0,\ldots,n+1$, are called {spectral bundles} and they are
independent of the choice of $\gamma_0,\ldots,\gamma_n$ in
{\rm\x{bc1}}, {\rm\x{bc2}} and {\rm\x{bc3}}. Moreover
\[\cW_0\oplus\cdots\oplus\cW_{n+1}=\Z \times \R^{N}\] is a Whitney
sum, i.e.,  $\cW_i \cap \cW_j=\Z \times \{0\}$ for $i\neq j$ and
$\cW_0+\cdots+\cW_{n+1}=\Z \times \R^{N}$.
\end{theorem}

\noindent\prf{Recall that the resolvent set $\rho_{NED}(A)$
is open and therefore $\Sigma_{NED}(A)$ is
the disjoint union of closed intervals. Next we will show that
$\Sigma_{NED}(A)$ consists of at most $N$ intervals. Indeed, if
$\Sigma_{NED}(A)$ contains $N+1$ components, then one can choose a
collections of points $\zeta_{1},\ldots,\zeta_{N}$ in
$\rho_{NED}(A)$ such that $\zeta_{1}< \cdots < \zeta_{N}$ and each
of the intervals $(0,\zeta_{1}),(\zeta_{1},\zeta_{2}),\ldots,$
$(\zeta_{N-1},\zeta_{N}),(\zeta_{N},\infty)$ has nonempty
intersection with the spectrum $\Sigma_{NED}(A)$. Now alternative II
of Lemma \ref{lemma2.3} implies \[0 \leq \dim \cS_{\zeta_1} < \cdots
< \dim \cS_{\zeta_N} \leq N\] and therefore either $\dim
\cS_{\zeta_1}=0$ or $\dim \cS_{\zeta_N}=N$ or both. Without loss of
generality, $\dim \cS_{\zeta_N}=N$, i.e., $\cS_{\zeta_N}=\Z \times
\R^{N}$. Assume that \[x_{k+1} =\tfrac{1}{\zeta_N}A_{k}x_{k}\] admits a strong
\ned~ with invariant projector $P \equiv {\rm Id}$, then \[x_{k+1}
=\tfrac{1}{\zeta}A_{k}x_{k}\] also admits a strong \ned~ with the same
projector for every $\zeta>\zeta_N$. Now we have the conclusion
$(\zeta_{N},\infty)\subset \rho_{NED}(A)$, which is a contradiction.
This proves the alternatives for $\Sigma_{NED}(A)$.

Due to Lemma \ref{lemma2.3}, the sets $\cW_0,\ldots,\cW_{n+1}$ are
invariant vector bundles. To prove now that $\dim \cW_1\geq
1,\ldots, \dim \cW_n\geq 1$ for $n \geq 1$, let us assume that $\dim
\cW_1=0$, i.e., $\cU_{\gamma_0} \cap \cS_{\gamma_1}=\Z \times
\{0\}$. If $(0,b_1]$ is a spectral interval this implies that
$\cS_{\gamma_1}=\Z \times \{0\}$. Then the projector of the \ned~ of
\[x_{k+1} =\tfrac{1}{\gamma_1}A_{k}x_{k}\] is $0$ and then we get
the contraction $(0,\gamma_1)\subset \rho_{NED}(A)$. If $[a_1,b_1]$
is a spectral interval then $[\gamma_0,\gamma_1]\cap \Sigma_{NED}(A)
\neq \emptyset$ and alternative II of Lemma \ref{lemma2.3} yields a
contradiction. Therefore $\dim \cW_1\geq 1$ and similarly $\dim
\cW_n\geq 1$. Furthermore for $n\geq 3$ and $i=2,\ldots,n-1$ one has
$[\gamma_{i-1},\gamma_{i}]\cap \Sigma_{NED}(A) \neq \emptyset$ and
again alternative II of Lemma \ref{lemma2.3} yields $\dim \cW_i\geq
1$.

For $i<j$ we have $\cW_{i}\subset \cS_{\gamma_i}$ and
$\cW_{i}\subset \cU_{\gamma_{j-1}} \subset \cU_{\gamma_i}$ and with
Lemma \ref{lem2.1} this gives $\cW_{i} \cap \cW_{j}\subset
\cS_{\gamma_i}\cap \cU_{\gamma_i}=\Z \times \{0\}$, so $\cW_{i} \cap
\cW_{j} =\Z \times \{0\}$ for $i \neq j$.

To show that $\cW_0\oplus\cdots\oplus\cW_{n+1}=\Z \times \R^{N}$,
recall the monotonicity relations $\cS_{\gamma_0}\subset \cdots
\subset \cS_{\gamma_n}$,
$\cU_{\gamma_0}\supset\cdots\supset\cU_{\gamma_n}$, and the identity
$\cS_\gamma \oplus \cU_\gamma = \Z \times \R^N$ for $\gamma \in
\R^{+}$. Therefore $\Z \times \R^{N}=\cW_0 \times \cU_{\gamma_0}$.
Now we have
 \beaa
\Z \times \R^{N} \EQ \cW_0 + \cU_{\gamma_0} \cap [\cS_{\gamma_1}+\cU_{\gamma_1}]\\
\EQ \cW_0 + [\cU_{\gamma_0} \cap \cS_{\gamma_1}]+\cU_{\gamma_1}\\
\EQ \cW_0 + \cW_1+\cU_{\gamma_1}. \eeaa
Doing the same for $\cU_{\gamma_1}$, we get
\beaa
\Z \times \R^{N} \EQ \cW_0 + \cW_1 + \cU_{\gamma_1} \cap [\cS_{\gamma_2}+\cU_{\gamma_2}]\\
\EQ \cW_0 + \cW_1 + [\cU_{\gamma_1} \cap \cS_{\gamma_2}]+\cU_{\gamma_2}\\
\EQ \cW_0 + \cW_1 + \cW_2+\cU_{\gamma_2}, \eeaa and mathematical
induction yields  $\Z \times \R^{N}=\cW_0+\cdots+\cW_{n+1}$. To
finish the proof, let $\tilde{\gamma}_0,\ldots,\tilde{\gamma}_n \in
\rho_{NED}(A)$ be given with the properties  {\rm\x{bc1}},
{\rm\x{bc2}} and {\rm\x{bc3}}. Then  alternative I of Lemma
\ref{lemma2.3} implies
\[\cS_{\gamma_i}=\cS_{\tilde{\gamma}_i} \quad {\rm and} \quad \cU_{\gamma_i}=\cU_{\tilde{\gamma}_i} \quad {\rm for} \quad i=0,\ldots,n\]
and therefore the invariant vector bundles $\cW_0,\ldots,\cW_{n+1}$
are independent of the choice of $\gamma_0,\ldots,\gamma_n$ in
{\rm\x{bc1}}, {\rm\x{bc2}} and {\rm\x{bc3}}.} $\Box$

\begin{definition} \label{Def2.2}
We say that \emph{\x{le}} is nonuniformly exponentially bounded
if there exist constants $K >0,\varepsilon\geq 1$ and $a \geq 1$
such that
\begin{equation}  \label{2.2}
 \|\Phi(k,l)\|
 \leq
Ka^{|k-l|}\varepsilon^{l},\quad k,l\in\Z.
\end{equation}
\end{definition}

\begin{lemma}\label{lemma2.4}
Assume that {\rm \x{le}} is nonuniformly exponentially bounded. Then
$\Sigma_{NED}(A)$ is a bounded closed set and
$\Sigma_{NED}(A)\subset [\frac{1}{a},a]$.
\end{lemma}

\noindent\prf{Assume that (\ref{2.2})
holds. Let $\gamma > a$ and $0<\alpha :=
\frac{a}{\gamma}<1$, then estimate (\ref{2.2}) implies
\[ \|\Phi_{\gamma}(k,l)\|
 \leq
K\alpha^{k-l}\varepsilon^{l}, \quad k\geq l.\]
Therefore \x{le} admits a \ned~ with invariant projector $P = I$. We
have $\gamma \in \rho_{NED}(A)$ and similarly for $0<\gamma <
\frac{1}{a}$, therefore $\Sigma_{NED}(A) \subset [\frac{1}{a},a]$.}
$\Box$

\begin{corollary}\label{BG}
 If {\rm\x{le}} is nonuniformly exponentially bounded, then the nonuniform dichotomy spectrum $\Sigma_{NED}(A)$
of {\rm\x{le}} is the disjoint union of $n$ closed intervals
 where $0 \leq n \leq N$, i.e.,
 \[\Sigma_{NED}(A) =
      [a_1,b_1]
    \cup [a_2,b_2] \cup \cdots \cup [a_{n-1},b_{n-1}]
     \cup
    [a_n,b_n],\]
where $a_1 \leq b_1 < a_2 \leq b_2 < \cdots < a_n \leq b_n$.
 \end{corollary}

From Proposition \ref{pro2.1}, we know $\Sigma_{NED}(A)\subset
\Sigma_{ED}(A)$. Finally in this Section, we present an example to illustrate that
$\Sigma_{NED}(A)\neq\Sigma_{ED}(A)$ can occur.

\begin{example}\lb{exa-2} Given $\omega>5a>0$. Consider the scalar
equation \be\lb{4.3}u_{k+1}=A_{k}u_{k}\ee with
\[A_{k}=e^{-\omega+ak(-1)^{k}-a(k-1)(-1)^{(k-1)}}.\]Then
$\Sigma_{NED}(A)= [e^{-\omega-a},e^{-\omega+a}]$ and
$\Sigma_{ED}(A)=\R^+$.
\end{example}

\noindent\prf{The
evolution operator of \x{4.3} is given by
\[\Phi(k,l)=e^{-\omega(k-l-1)-a(k-l-1)(-1)^{k-1}-al(-1)^{(k-1)}+al(-1)^{l}}.\]
For any $\gamma\in \R^+$ the evolution operator of the equation
\be\lb{uk}u_{k+1}=\tfrac{1}{\gamma}A_{k}u_{k}\ee  is given by
\begin{equation}\label{4.4}
\Phi_{\gamma}(k,l)=(\tfrac{1}{\gamma})^{(k-l)}e^{-\omega(k-l-1)-a(k-l-1)(-1)^{k-1}-al(-1)^{(k-1)}+al(-1)^{l}}.
\end{equation}
For any $\gamma\in (e^{(-\omega+5a)},+\infty)$,  it follows from
(\ref{4.4}) that \be\lb{4.5}|\Phi_{\gamma}(k,l)| \leq
e^{\omega-a}\bigg{(}\frac{e^{-\omega+a}}{\gamma}\bigg{)}^{k-l
}e^{2al},\,\,k\geq l,\ee which implies that the equation \x{uk}
admits a strong \ned~ with $P={\rm Id}$, by taking
\[K=e^{\omega-a},~~\alpha=\frac{e^{-\omega+a}}{\gamma}<1,~~\varepsilon=e^{2a}>0.\]
Thus,
 \begin{equation}\label{4.6}
(e^{-\omega+5a},+\infty)\subset \rho_{NED}(A).
 \end{equation}
For any $\widetilde{\gamma}\in (0, e^{-\omega-5a})$,  it follows from (\ref{4.4}) that
 \be\lb{4.7}
|\Phi_{\gamma}(k,l)| \leq
e^{\omega+a}\bigg{(}\frac{e^{-\omega-a}}{\gamma}\bigg{)}^{k-l
}e^{2al},\,\,k\leq l,\ee which implies that \x{uk} admits a strong \ned~ with $P=0$,
  by taking \[K=e^{\omega+a},~~\alpha=\frac{\gamma}{e^{-\omega-a}}<1,~~\varepsilon=e^{2a}>0.\]
 Thus,
 \begin{equation}\label{4.8}
(0, e^{-\omega-5a})\subset \rho_{NED}(A).
  \end{equation}
 It follows from (\ref{4.6}) and (\ref{4.8}) that
 \[(0, e^{-\omega-5a}) \cup   (e^{-\omega+5a},+\infty)\subset
\rho_{NED}(A),\]
 which implies that
 \[\Sigma_{NED}(A) \subset[e^{-\omega-5a},e^{-\omega+5a}].\]
Next we show that
 \[[e^{-\omega-5a},e^{-\omega+5a}]\subset\Sigma_{NED}(A).\]
To do this, we first prove that $\gamma_{1}=e^{-\omega+5a} \in
 \Sigma_{NED}(A)$. The evolution operator of the system
\[u_{k+1}=\tfrac{1}{\gamma_{1}}A_{k}u_{k}\]  is given as
\[\Phi_{\gamma_{1}}(k,l)=e^{\omega-a}e^{-a(k-l-1)(1+(-1)^{k-1})-al(-1)^{(k-1)}+al(-1)^{l}}.\]
It is easy to see that there do not exist $K$, $\alpha>0$ and
$\varepsilon>0$ such that
\[\|\Phi_{\gamma_{1}}(k,l)\| \leq
 K\alpha^{k-l}\varepsilon^{l}, \ \ \ \ \mbox{ for} \,\,\, k\geq l,\]
 or
\[\|\Phi_{\gamma_{1}}(k,l)\|
 \leq
K(\tfrac{1}{\alpha})^{k-l}\varepsilon^{l},\ \ \mbox{ for} \,\,\,
k\leq l.\]

Therefore $\gamma_{1}=e^{-\omega+5a} \in
 \Sigma_{NED}(A)$. In a similar manner, we can prove $\gamma_{2}=e^{-\omega-5a} \in
 \Sigma_{NED}(A)$. We can see from Theorem \ref{main21} that \x{4.3}
 has at most one nonuniform dichotomy spectral interval, which means that $
[e^{-\omega-5a},e^{-\omega+5a}]\subset\Sigma_{NED}(A)$ and therefore
$[e^{-\omega-5a},e^{-\omega+5a}]=\Sigma_{NED}(A).$

On the other hand, using a similar argument as in equations \x{4.1},
we know that the nonuniform part  $\varepsilon^{l}$ cannot be removed in
the estimates \x{4.5} and \x{4.7}. Therefore, \x{4.3} does not admit an exponential
dichotomy, which means that $\Sigma_{ED}(A)=\R^+.$ } $\Box$

\section{Reducibility}\setcounter{equation}{0}

In this section we employ Theorem \ref{main21} to prove a reducibility result. For the
reducibility results in the setting of an exponential dichotomy, we refer the reader to \cite{cop67, pal,
s-jlms} and the references therein.

\begin{lemma}\label{lemma3.2}  The projector of equation {\em \x{le}}
can be chosen as $\tilde{P}= \left(\begin{array}{lll} I_{N_1} & 0
\\
0 &  0_{N_2}
\end{array}
\right)$ with $N_1 = {\rm dim\, im} \tilde{P}$ and $N_2 = {\rm dim\,
ker} \tilde{P}$, and the fundamental matrix $X_{k}$ can be chosen
suitably such that the estimates  {\em\x{npl}-\x{nql}}  can be
rewritten as \be\lb{npl-e}\|X_{k}\tilde{P}X^{-1}_{l}\|\leq
K\alpha^{k-l}\varepsilon^{l}, \quad k\geq l,\ee and
\be\lb{nql-e}\|X_{k}\tilde{Q}X^{-1}_{l}\|\leq
K(\frac{1}{\alpha})^{k-l}\varepsilon^{l},\quad k\leq l,\ee where
$\tilde{Q}={\rm Id}-\tilde{P}.$
\end{lemma}

\noindent\prf{Let $n\in \mathbb{Z}$ be arbitrary but fixed. Note
that the rank of the projector $P_n$ is independent of $n\in
\mathbb{Z}$ (see \cite[Page 1100]{bds}), then there exists a
nondegenerate matrix $T\in \mathbb{R}^{N\times N}$ such that
\[\tilde{P}:=\left(\begin{array}{lll} I_{N_1} & 0
 \\
 0 &  0_{N_2}
 \end{array}
 \right)=TP_{n}T^{-1} \]
 with $N_1 = {\rm dim\, im} \tilde{P}$ and $N_2 = {\rm dim\,
ker} \tilde{P}$. Define
\[X_{k} := \Phi(k,n)T^{-1} \quad \text{for } k \in \mathbb{Z}
\quad \text{and} \quad
  \tilde{P}:=  \left(\begin{array}{lll}
I_{N_1} & 0
 \\
 0 &  0_{N_2}
 \end{array}
 \right)= T P_{n} T^{-1}.\]
Then
 \begin{equation}
\|X_{k}\tilde{P}X^{-1}_{l}\|=\|\Phi(k,n)T^{-1} \tilde{P}T
\Phi^{-1}(l,n)\|=\|\Phi(k,n)P_{n}\Phi^{-1}(l,n)\|. \label{3.7}
\end{equation}
On the other hand, we have
 \bea\lb{EQ}
\|\Phi(k,l)P_{l}\|\EQ\|\Phi(k,n)\Phi(n,l)P_{l}\|\nn\\
\EQ\|\Phi(k,n)P_{n}\Phi(n,l)\|\nn\\
\EQ\|\Phi(k,n)P_{n}\Phi^{-1}(l,n)\|. \eea It follows from
(\ref{3.7}) and (\ref{EQ}) that \x{npl}-\x{nql} can be rewritten in
the form \x{npl-e}-\x{nql-e}.} $\Box$

Now we recall the definition of kinematic similarity and
several results in Coppel \cite{cop} and Aulbach et al. \cite{amz}.
\begin{definition}
Equation {\rm(\ref{le})} is said to be kinematically similar to
another equation \be \lb{3.1}y_{k+1}=B_{k}y_{k} \ee with $k\in \Z$,
if there exists an  invertible matrix $S_{k}$ with $\|S_{k}\|\leq M$
and $\|S_{k}^{-1}\|\leq M(M>0),$ which satisfies the difference
equation
\[S_{k+1}B_{k}=A_{k}S_{k}.\]
The change of variables $x_{k}=S_{k}y_{k}$ then transforms {\rm(\ref{le})} into {\rm(\ref{3.1})}.
\end{definition}

The next lemma is important to establish the reducibility results
and its proof follows along the lines of the proof of Siegmund
\cite{s-jlms}. See also Coppel \cite{cop} and Aulbach et al.
\cite{amz}

\begin{lemma}\label{lemma 3.1}{\rm \cite[Chapter 5]{cop}} Let $P$
be an orthogonal projection ($P^{T}=P$) and let $X$ be an invertible matrix.
 Then there exists an invertible matrix function $S:\mathbb{Z}\rightarrow\R^{N\times N}$ such that
\[S_{k}PS^{-1}_{k}=X_{k}PX^{-1}_{k},\,\,\,\,\,\,\, \,\,\,\,\,\,\,
S_{k}QS^{-1}_{k}=X_{k}QX^{-1}_{k},\] and
\begin{equation*}
\|S_{k}\|{\leq}\sqrt{2},\label{3.2}
\end{equation*}
\begin{equation*}
\|S^{-1}_{k}\|{\leq}\big[\|X_{k}PX^{-1}_{k}\|^{2}+\|X_{k}(I-P)X^{-1}_{k}\|^{2}\big]^{\frac{1}{2}}\label{3.3},
\end{equation*} where $k\in \Z$ and $Q={\rm Id}-P.$
Define\[\widetilde{R} :\mathbb{Z}\rightarrow \mathbb{R}^{N \times
N},
      \;
k \mapsto PX_{k}^T X_{k} P+ [{\rm Id} - P] X_{k}^T X_{k} [{\rm Id} -
P].\] Then the mapping is a positive definite, symmetric
 matrix for every $k \in \mathbb{Z}$. Moreover there is a unique function
 \[R : \mathbb{Z}\rightarrow \mathbb{R}^{N \times N}\]
 of positive definite symmetric matrices $R_{k}$, $k \in \mathbb{Z}$, with
 \[R_{k}^2 = \widetilde{R}_{k}
 ,
 \quad
 PR_{k} = R_{k} P\;.\]
\end{lemma}

We remark that $S^{-1}_{k}$ in Lemma \ref{lemma 3.1} is bounded
in the setting of an exponential dichotomy. However, in the
setting of a \ned, $S^{-1}_{k}$ can be unbounded, because
$\|\Phi(k,k)P_{k}\|\leq K\varepsilon^{k}$ for $k\geq 0$. To overcome
the difficulty, we introduce a new version of non-degeneracy,
so-called weak non-degeneracy and define the concept of weak
kinematical similarity. Some results will be obtained on the
decoupling into two blocks which will play an important role in the
analysis of reducibility.

\begin{definition} \label{Def3.2}
$S:\mathbb{Z}\rightarrow\R^{N\times N}$ is called  {\rm weakly
non-degenerate}  if  there exists a constant $M=M(\varepsilon)>0$ such
that
\[\|S_{k}\|\leq M\varepsilon^{|k|}\,\,\,\,\mbox{and}\,\,\,\,  \|S^{-1}_{k}\|  \leq M\varepsilon^{|k|},
 \,\,\,\,\mbox{for all}\,\,\,\, k\in \mathbb{Z}.\] \end{definition}

\begin{definition} \label{Def3.3}
If there exists a weakly non-degenerate matrix $S_{k}$ such that
\[S_{k+1}B_{k}=A_{k}S_{k},\]
then equation {\rm(\ref{le})} is weakly kinematically similar to
equation {\rm(\ref{3.1})}.
 For short, we denote {\rm(\ref{le})} $\overset{w}{\sim}$ {\rm(\ref{3.1})} or $A_{k} \overset{w}{\sim}  B_{k}$.
\end{definition}
For the sake of comparison, we denote kinematical similarity by
{\rm(\ref{le})} $ \sim $ {\rm(\ref{3.1})} or $A_{k} \sim B_{k}$.

\begin{definition} \label{Def3.4}
 We say that equation {\rm(\ref{le})} is reducible, if it is weakly kinematically similar
 to equation {\rm(\ref{3.1})}
 whose coefficient matrix $B_{k}$ has the block form
 \begin{equation}
 \left(\begin{array}{lll}
 B^1_{k} & 0
 \\
 0 &  B^2_{k}
 \end{array}
 \right),
 \label{3.4}
 \end{equation}
 where $B^1_{k}$ and $B^2_{k}$ are matrices of smaller size than $B_{k}$.
\end{definition}

The following theorem shows that if (\ref{le}) admits a nonuniform
exponential dichotomy, then there exists a weakly non-degenerate
transformation such that $A_{k}\overset{w}{\sim} B_{k}$ and $B_{k}$
has the block form (\ref{3.4}),
 i.e., system (\ref{le}) is reducible.

\begin{theorem}\label{Th3.1}
Assume that {\rm(\ref{le})} admits a nonuniform exponential dichotomy {\rm(}not necessary strong{\rm)}
of the form {\rm\x{npl-e}-\x{nql-e}} with invariant projector
$P_k\neq0,{\rm Id}$. Then {\rm(\ref{le})} is weakly
kinematically similar to a decoupled system
\begin{equation}  \label{3.5}
x_{k+1} =
\begin{pmatrix}
B^1_{k} & 0
 \\
 0 &  B^2_{k}
\end{pmatrix}
x_{k}
\end{equation}
for some locally integrable matrix functions
\[B^1 :  \mathbb{Z} \rightarrow  \mathbb{R}^{N_1 \times N_1}
\quad \text{and} \quad B^2 :  \mathbb{Z} \rightarrow \mathbb{R}^{N_2
\times N_2}\] where $N_1 := \dim \im \tilde{P}$ and $N_2 := \dim
\ker \tilde{P}$. That is, system {\rm(\ref{le})} is reducible.
\end{theorem}

\noindent\prf{Since equation (\ref{le}) admits a nonuniform
exponential dichotomy of the form \x{npl}-\x{nql} with invariant
projector $P_{k}\neq0,{\rm Id}$, by Lemma \ref{lemma3.2}, we can choose
suitable fundamental matrix $X_{k}$
  and the projector $\tilde{P}= \left(\begin{array}{lll}
I_{N_1} & 0
 \\
 0 &  0
 \end{array}\right)$,$(0<N_{1}<N)$ such that the estimates \x{npl-e}-\x{nql-e} hold.
By Lemma \ref{lemma 3.1} and the estimates \x{npl-e}-\x{nql-e},
there exists a $M=M(\varepsilon)>0$ large enough such that
\[\|S_{k}\|{\leq}\sqrt{2}\leq M\varepsilon^{|k|},\]
\[\|S^{-1}_{k}\|{\leq}\big[\|X_{k}\tilde{P}X^{-1}_{k}\|^{2}+\|X_{k}(I-\tilde{P})X^{-1}_{k}\|^{2}\big]^{\frac{1}{2}}\leq
\sqrt{2}K\varepsilon^{|k|}.\]
 Thus, $S$ is weakly
non-degenerate. Setting
 \[B_{k}=R_{k+1}R^{-1}_{k},\]
 where $R_{k}$ is defined in Lemma \ref{lemma 3.1} and $X_{k}=S_{k}R_{k}$.
 Obviously, $R_{k}$ is the fundamental matrix of linear system
\begin{equation*}
y_{k+1}=B_{k}y_{k}. \label{2.11}
\end{equation*}
Now we need to show that $A_{k}\overset{w}{\sim}  B_{k}$ and $B_{k}$
has the block diagonal form
\begin{equation*}
B_{k}= \left(\begin{array}{lll} B^1_{k} & 0
 \\
 0 &  B^{2}_{k}
 \end{array}
 \right),\,\,\,\,\mbox{for} \,\,\,\, k\in \mathbb{Z}.
\label{2.12}
\end{equation*}
 First, we show that $A_{k}\overset{w}{\sim} B_{k}$. In fact,
\beaa
S_{k+1}B_{k}\EQ X_{k+1}R^{-1}_{k+1}B_{k}\\
\EQ A_{k}X_{k}R^{-1}_{k}B^{-1}_{k}B_{k}\\
\EQ A_{k}S_{k},\eeaa
which implies that $A_{k}\overset{w}{\sim}  B_{k}$.

Now we show that system {\rm(\ref{le})} is weakly kinematically
similar to (\ref{3.5}). By Lemma \ref{lemma 3.1}, $R_{k+1}$ and
$R_{k}^{-1}$ commute with the matrix $\tilde{P}$ for every $k \in
\mathbb{Z}$. It follows that \be\lb{2.13} \tilde{P} B_{k} = B_{k}
\tilde{P}\ee for all $k \in \mathbb{Z}$. Now we decompose
$B_{k}: \Z \rightarrow \R^{N \times N}$ into four functions
 \[\begin{array}{lll}
B^{1}_{k}: \Z \rightarrow \R^{N_{1} \times N_{1}}, & B^{2}_{k}: \Z
\rightarrow \R^{N_{2} \times N_{2}},\\
 B^{3}_{k}: \Z \rightarrow \R^{N_{1} \times N_{2}}, & B^{4}_{k}: \Z
\rightarrow \R^{N_{2} \times N_{1}},
   \end{array}\]
 with
\[B_{k}= \left(\begin{array}{lll} B^1_{k} & B^3_{k}
 \\
B^4_{k} &  B^{2}_{k}
 \end{array}
 \right), \quad \,\,\,\, k\in \mathbb{Z}.\]
Identity \x{2.13} implies that
\[\left(\begin{array}{lll} B^1_{k} & B^3_{k}
 \\
0 &  0
 \end{array}
 \right)=\left(\begin{array}{lll} B^1_{k} & 0
 \\
B^4_{k} &  0
 \end{array} \right), \quad \,\,\,\, k\in \mathbb{Z}.\]
Therefore $B^3_{k}\equiv 0$ and $B^4_{k}\equiv 0$. Thus $B_{k}$ has
the block form
\[B_{k}= \left(\begin{array}{lll} B^1_{k} & 0
 \\
0 &  B^{2}_{k}
 \end{array} \right), \quad \,\,\,\, k\in \mathbb{Z}.\]
Now the proof is finished.} $\Box$

From Theorem \ref{Th3.1}, we know that if (\ref{le}) admits a
nonuniform exponential dichotomy, then there exists a weakly
non-degenerate transformation $S_{k}$ such that
$A_{k}\overset{w}{\sim} B_{k}$ and $B_{k}$ has two blocks of the
form (\ref{3.4}).

\begin{lemma}\label{lemma3.3}  Assume that \emph{(\ref{le})} admits a nonuniform exponential
dichotomy with the form of estimates \emph{\x{npl-e}-\x{nql-e}} and
$rank(\tilde{P})=N_{1}, (0 < N_{1}< N)$, and there exists a weakly
non-degenerate transformation $S_{k}$ such that
$A_{k}\overset{w}{\sim} B_{k}$. Then system \emph{(\ref{3.1})}
also admits a nonuniform exponential dichotomy, and the projector
has the same rank.
\end{lemma}

\noindent\prf{Suppose that $S_{k}$ is weakly non-degenerate, which
means that  there exists $M=M(\varepsilon)>0$ such that $\|S_{k}\|\leq
M\varepsilon^{|k|}$ and $\|S^{-1}_{k}\|  \leq
M\varepsilon^{|k|}$ and such that $A_{k}\overset{w}{\sim}
B_{k}$. Let $X_{k}=S_{k}Y_{k}$. It is easy to see that $Y_{k}$ is
the fundamental matrix of system (\ref{3.1}). To prove that system
(\ref{3.1}) admits a nonuniform exponential dichotomy, we first
consider the case $k\geq l$ and obtain
\begin{equation}
\begin{array}{lll}
\|Y_{k}\tilde{P}Y^{-1}_{l}\|&=&
\|S^{-1}_{k}X_{k}\tilde{P}X^{-1}_{l}S_{k}\|
\\
&\leq&
 \|S^{-1}_{k}\|\cdot\|X_{k}\tilde{P}X^{-1}_{l}\|\cdot\|S_{l}\|
 \\
&\leq&
 K M^{2}\varepsilon^{|k|} \alpha^{k-l}\varepsilon^{l}\varepsilon^{|l|}
  \\
&\leq&
 K M_1^{2}(\varepsilon\alpha)^{k-l}\varepsilon^{l},\,\,k\geq l,
 \end{array}
\label{3.12}
\end{equation}where $M_1=M\varepsilon^{2|l|}.$
Similar argument shows that
 \begin{equation}
\|Y_{k}\tilde{Q}Y^{-1}_{l}\|\leq  K
M^{2}_{1}(\tfrac{1}{\varepsilon\alpha})^{k-l}\varepsilon^{l},\,\,k\leq l.\label{3.13}
\end{equation}
Form (\ref{3.12}) and (\ref{3.13}), it is easy to see that system
(\ref{3.1}) admits a nonuniform exponential dichotomy.  Clearly, the rank of the projector is
$k$.} $\Box$

\begin{lemma}\label{lemma3.4}  Assume that the systems \emph{(\ref{le})}
 and \emph{(\ref{3.1})} are weakly kinematically similar via $S_{k}$. If
 for a constant $\gamma \in \R^+$ the system {\rm\x{2.1}}
admits a strong nonuniform exponential dichotomy with constants $K>0$,
$0<\alpha<1$, $\varepsilon \geq 1$ and invariant projector
$P$, then the system\be\lb{ykk}y_{k+1}=\tfrac{1}{\gamma}B_{k}y_{k}\ee also admits a strong
nonuniform exponential dichotomy.
 \end{lemma}

\noindent\prf{Obviously, $P$ is also an invariant
projector for \x{le}. The dichotomy estimates are equivalent to
\[\|X_{k}P X^{-1}_{l}\| \leq K\alpha^{k-l}\varepsilon^{l}, \quad k\geq l,\]
and \[\|X_{k}P X^{-1}_{l}\|\leq
K(\tfrac{1}{\alpha})^{k-l}\varepsilon^{l},\quad k\leq l.\] Using
Lemma \ref{lemma3.3}, it is easy to see that
\[\|Y_{k}P Y^{-1}_{l}\| \leq
 K'_\gamma(\varepsilon\alpha)^{k-l}\varepsilon^{l}, \quad k\geq l,\]
and\[\|Y_{k}P Y^{-1}_{l}\|\leq
K'_\gamma(\tfrac{1}{\varepsilon\alpha})^{k-l}\varepsilon^{l},\quad k\leq
l,\] for some constant $K'_\gamma\geq 1.$ Therefore, \x{ykk} admits a strong nonuniform exponential dichotomy.} $\Box$

The following result follows directly from Lemma \ref{lemma3.4}.
\begin{corollary}\label{cor3.1}  Assume that there exists a weakly non-degenerate transformation $S_{k}$
such that $A_{k}\overset{w}{\sim} B_{k}$. Then $\Sigma_{NED}(A) =
  \Sigma_{NED}(B)$,
 i.e.,\[\Sigma_{NED}(A) =
  \left\{ \begin{matrix}
    [a_1,b_1]
    \\
    \text{or}
    \\
    (0,b_1]
  \end{matrix}\right\}
  \cup [a_2,b_2] \cup \cdots \cup [a_{n-1},b_{n-1}] \cup
  \left\{\begin{matrix}
    [a_n,b_n]
    \\
    \text{or}
    \\
    {[}a_n,\infty)
  \end{matrix}\right\}=\Sigma_{NED}(B).\]
\end{corollary}

Now we are in a position to prove the reducibility result.
\begin{theorem}[Reducibility Theorem] Assume that {\rm(\ref{le})} admits a strong nonuniform
exponential dichotomy. Due to Theorem {\rm
\ref{main21}}, the dichotomy spectrum is either empty or the
disjoint union of $n$ closed spectral intervals $\cI_1, \dots,
\cI_n$ with $1 \leq n \leq N$, i.e.,
\[\Sigma_{NED}(A) = \emptyset \quad(n = 0)
  \qquad\text{ or }\qquad
  \Sigma_{NED}(A) = \cI_1 \cup \cdots \cup \cI_n\;.\]
Then there exists a weakly kinematic
similarity action $S : \Z \to \R^{N \times N}$ between
{\rm(\ref{le})} and a block diagonal system
\[  \label{nsystem} x_{k+1} =\begin{pmatrix}
    B^{0}_{k} & & \\
    & \ddots & \\
    & & B^{n+1}_{k}
  \end{pmatrix}
x_{k}\]with $B^{i} : \Z \rightarrow
\R^{N_i \times N_i}$, $N_i=\dim \cW_i$, and
\[\Sigma_{NED}(B^0) = \emptyset
  \;,\Sigma_{NED}(B^1) = \cI_1
  \;,
  \dots
  \;,
  \Sigma_{NED}(B^n) =\cI_n,\Sigma_{NED}(B^{n+1}) = \emptyset.\]
\end{theorem}
\noindent\prf{If for any $\gamma\in\R^+$,  system \x{2.1} admits a
strong \ned, then $\Sigma_{NED}(A) = \emptyset$. Conversely, for any
$\gamma\in \R^+$, system \x{2.1} does not admit a strong \ned,
then $\Sigma_{NED}(A) =\R^+$. Now, we prove the theorem for the
nontrivial case ($\Sigma_{NED}(A) \neq \emptyset$ and
$\Sigma_{NED}(A) \neq\R^+$).

Recall that the resolvent set
$\rho_{NED}(A)$ is open and therefore the
dichotomy spectrum $\Sigma_{NED}(A)$ is the disjoint union of closed
intervals. Using Theorem \ref{main21}, we can assume\[\cI_{1} =
\left\{\begin{matrix}
    [a_1,b_1]
    \\
    \text{or}
    \\
    (0,b_1]\end{matrix}\right\}
, \cI_{2}=[a_2,b_2], \ldots , \cI_{n-1}=[a_{n-1},b_{n-1}] ,
  \cI_{n} =\left\{\begin{matrix}
    [a_n,b_n]
    \\
    \text{or}
    \\
    {[}a_n,\infty)\end{matrix}\right\}\]
with $0<a_1 \leq b_1 < a_2 \leq b_2 < \ldots < a_n \leq b_n$.

If $\cI_{1}=[a_1,b_1]$ is a spectral interval, then $(0,
\gamma_0)\subset \rho_{NED}(A)$ and $\cW_0=\cS_{\gamma_0}$ for some
$\gamma_0 < a_1$ due to Theorem \ref{main21}, which implies that
\[x_{k+1}=\tfrac{1}{\gamma_0}A_{k}x_{k}\]
admits a strong \ned~ with an invariant projector $\tilde{P}_0$. By Theorem
\ref{Th3.1} and Corollary \ref{cor3.1}, there exists a weakly
non-degenerate transformation $x_k=S_k^{0}x_k^{(0)}$ with
$\|S_{k}^{0}\|\leq M_0\varepsilon^{|k|}$ and
$\|(S_{k}^{0})^{-1}\| \leq M_0\varepsilon^{|k|}$ for some positive constant $M_0=M_0(\varepsilon)$ and such
that $A_{k}\overset{w}{\sim} A^{0}_{k}$ and $A^{0}_{k}$ has two
blocks of the form $A^{0}_{k}=\left(\begin{array}{lll}
 B^{0}_{k} & 0
 \\
 0 &   B^{0,*}_{k}
 \end{array}\right)
$ with $\dim B_{k}^{0}=\dim \im \tilde{P}_0=\dim \cS_{\gamma_0}=\dim
\cW_{0}=:N_{0}$ due to Theorem~\ref{Th3.1}, Lemma~\ref{lem2.1} and
Theorem \ref{main21}. If $\cI_{1}=(0,b_1]$ is a spectral interval, a
block $B_{k}^{0}$ is omitted.

Now we consider the following system
\[x_{k+1}^{(0)}=A^{0}_{k}x_{k}^{(0)}=\left(\begin{array}{lll}
 B^{0}_{k} & 0
 \\
 0 &   B^{0,*}_{k}
 \end{array}\right)x_{k}^{(0)}.\]
By using Lemma \ref{lemma2.3}, we take $\gamma_1\in (b_{1},a_{2})$.
In view of $(b_{1},a_{2})\subset \rho_{NED}(B_k^{0,*})$,
$\gamma_1\in \rho_{NED}(B_k^{0,*})$, which implies that
\[x_{k+1}^{(0)}=\frac{1}{\gamma_1}\left(\begin{array}{lll}
 B^{0}_{k} & 0
 \\
 0 &   B^{0,*}_{k}
 \end{array}\right)x_{k}^{(0)}\]
admits a \ned ~with an invariant projector $\tilde{P}_1$. From the claim
above, we know that $\tilde{P}_1\neq 0,\,I$. Similarly by Theorem
\ref{Th3.1} and Corollary \ref{cor3.1}, there exists a weakly
non-degenerate transformation
\[x_k^{(0)}=S_k^{1}x_k^{(1)}=\left(\begin{array}{lll}
 I_{N_{0}} & 0
 \\
 0 & \tilde{S}^{1}_{k}
 \end{array}
 \right)x_k^{(1)}\]
with $\|\tilde{S}^{1}_{k}\|\leq M_1\varepsilon^{|k|}$
and $\|(\tilde{S}^{1}_{k})^{-1}\| \leq
M_1\varepsilon^{|k|}$ for some positive constant $M_1=M_1(\varepsilon)$ and such that
$B^{0,*}_{k}\overset{w}{\sim} \tilde{B}^{0,*}_{k}$ and $
\tilde{B}^{0,*}_{k}$ has two blocks of the form $
\tilde{B}^{0,*}_{k}=\left(\begin{array}{lll}
 B^{1}_{k} & 0
 \\
 0 &   B^{1,*}_{k}
 \end{array}
 \right)
$ with $\dim B_{k}^{1}=\dim \im \tilde{P}_1=\dim
\cS_{\gamma_1}\geq\dim(\cU_{\gamma_0}\cap \cS_{\gamma_1})=\dim
\cW_{1}=:N_{1}$ due to Theorem~\ref{Th3.1}, Lemma~\ref{lem2.1} and
Theorem \ref{main21}. In addition, using Theorem \ref{Th3.1} and Corollary \ref{cor3.1}, we have
\[\Sigma_{NED}(B_k^1) =\left\{\begin{matrix}
    [a_1,b_1]
    \\
    \text{or}
    \\
    (0,b_1]
  \end{matrix}\right\},~~
\Sigma_{NED}(B_k^{1,*}) =[a_2,b_2] \cup \cdots \cup
[a_{n-1},b_{n-1}] \cup
  \left\{\begin{matrix}
    [a_n,b_n]
    \\
    \text{or}
    \\
    {[}a_n,\infty)
  \end{matrix}\right\}.\]Now we can construct a weakly
non-degenerate transformation $x_k=\tilde{S}_kx_k^{(1)}$ with $
 \tilde{S}_{k}=S^{0}_{k}S^{1}_{k}=S^{0}_{k}\left(\begin{array}{lll}
 I_{N_{0}} & 0
 \\
 0 & \tilde{S}^{1}_{k}
 \end{array}\right)
$, where $\| \tilde{S}_{k}\|\leq M_0
M_1\varepsilon^{2|k|}$ and $\|
\tilde{S}_{k}^{-1}\| \leq M_0
M_1\varepsilon^{2|k|}$. Then
$A_{k}\overset{w}{\sim} A^{1}_{k}$ and $A^{1}_{k}$ has three blocks
of the form \[ A^{1}_{k}= \begin{pmatrix}
    B^{0}_{k} & & \\
    & B^{1}_{k} & \\
    & & B^{1, *}_{k}
\end{pmatrix}.\]

Applying similar procedures to $\gamma_{2}\in(b_{2},a_{3})$,
$\gamma_{3}\in(b_{3},a_{4}),\ldots$, we can construct a weakly
non-degenerate transformation $x_k=S_kx_k^{(n+1)}$ with \[
 S_{k}=S^{0}_{k}\left(\begin{array}{lll}
 I_{N_{0}} & 0
 \\
 0 & \tilde{S}^{1}_{k}
 \end{array}
 \right)\left(\begin{array}{lll}
 I_{N_{0}+N_{1}} & 0
 \\
 0 & \tilde{S}^{2}_{k}
 \end{array}
 \right)\cdots \left(\begin{array}{lll}
 I_{N_{0}+\ldots+N_{n-1}} & 0
 \\
 0 & \tilde{S}^{n}_{k}
 \end{array}
 \right)\]
such that $\| S_{k}\|\leq M_{\varepsilon} \varepsilon^{n|k|}$ and $\|
S_{k}^{-1}\| \leq M_{\varepsilon} \varepsilon^{n|k|}$ with
$M_{\varepsilon}=M_0\times \cdots \times
M_n$. Now
we can prove
\[A_{k}\overset{w}{\sim} A^{n}_{k}:=B_{k} =\begin{pmatrix}
B^{0}_{k} & & \\
& \ddots & \\
& & B^{n+1}_{k}\end{pmatrix}\]
with locally integrable functions $B^{i} : \Z \rightarrow \R^{N_i
\times N_i}$ and\[\Sigma_{NED}(B^0) = \emptyset
\;,\Sigma_{NED}(B^1) = \cI_1
\;,
\dots
\;,
\Sigma_{NED}(B^n) =\cI_n,\Sigma_{NED}(B^{n+1}) = \emptyset.\]

Finally, we show that $N_i=\dim \cW_i$. From the claim above, we note
that $\dim B^{0}_k=\dim \cW_0, ~\dim B^{1}_k\geq \dim \cW_1,\ldots,
\dim B^{n}_k\geq \dim \cW_n, ~\dim B^{n+1}_k=\dim \cW_{n+1}$ and
with Theorem \ref{main21} this gives $\dim \cW_0 + \cdots +\dim \cW_{n+1}=N$,
so $\dim B^{i}_k=\dim \cW_i$ for $i=0,\ldots,n+1$. Now the proof is finished.} $\Box$

\end{document}